%% file: rr.tex
\documentclass[11pt]{article}

\usepackage{amsmath,amssymb,amscd}
\usepackage[dvips]{graphics}

\newtheorem{theorem}{Theorem}[section]

\newtheorem{proposition}[theorem]{Proposition}

\newtheorem{lemma}[theorem]{Lemma}
\newtheorem{example}[theorem]{Example}

\newtheorem{remark}[theorem]{Remark}

\newtheorem{corollary}[theorem]{Corollary}

\newcommand{\ethrm}{\hspace*{\fill}
                      $\Box$
                      \vspace{.1in}}

\newcommand{\bp}{{\mathbb P}}

\newcommand{\bq}{{\mathbb Q}}

\newcommand{\ba}{{\mathbb A}}

\renewcommand{\L}{{\mathcal L}}
\renewcommand{\P}{{\mathcal P}}
\newcommand{\Q}{{\mathcal Q}}
\newcommand{\D}{{\mathcal D}}
\newcommand{\bk}{{\mathcal K}}
\newcommand{\bl}{{\mathcal L}} 
\newcommand{\bt}{{\mathcal T}}
\newcommand{\bn}{{\mathcal N}}

\newcommand{\bigo}{{\mathcal O}}

\newcommand{\ol}{\overline}
\newcommand{\ep}{\varepsilon}

\newcommand{\cy}{Calabi--Yau }

\newcommand{\lcm}{\mathop{\rm lcm}\nolimits}
\newcommand{\hcf}{\mathop{\rm hcf}\nolimits}

\newcommand{\proj}{\mathop{\rm Proj\,}\nolimits}

\newcommand{\sing}{\mathop{\rm Sing}\nolimits}

\newcommand{\wgr}{\mathop{\rm wGr}\nolimits}

\newcommand{\wogr}{\mathop{\rm wOGr}\nolimits}
\newcommand{\pf}{\mathop{\rm Pf}\nolimits}

\renewcommand{\Box}{\square}

\newcommand{\color}[6]{}

\begin{document}

\begin{center}
{\LARGE Orbifold Riemann--Roch for threefolds} \\
\vspace{0.15in}
{\LARGE with an application to Calabi--Yau geometry} \\
\vspace{0.2in}
{\large Anita Buckley and Bal\'{a}zs Szendr\H{o}i}\\
\vspace{0.15in}
{\large June 2004}
\vspace{0.1in}
\end{center}

{\small
\begin{center} {\sc abstract} \end{center}
{\leftskip=30pt \rightskip=30pt
We prove an orbifold Riemann--Roch formula for a polarized 3--fold $(X,D)$. 
As an application, we construct new families of projective Calabi--Yau 
threefolds.
\par}}

\section*{Introduction}
\label{introduction}
The aim of this paper is to state and prove a Riemann--Roch (RR) formula for a 
pair $(X,D)$ of a normal projective threefold $X$ with quotient singularities and
a $\bq$-Cartier Weil divisor $D$. Under certain conditions, 
we prove the existence of a formula
$$\chi(X,\bigo_X(D))=\mbox{ RR--type expression in }D+\sum_P c_P(D)+\sum_C s_C(D),$$
where $c_P(D)$ is the contribution from a singular point $P$ of $(X,D)$ 
and $s_C(D)$ is the contribution from a curve $C$ of singularities. 
These contributions depend on the type of the quotient singularities and on the 
embedding of $\sing X$ in $X$.

Explicit formulae of this type for surfaces and threefolds with isolated 
canonical singularities were first studied by Fletcher~\cite{fletcher1} and 
Reid~\cite{ypg}. Earlier Kawasaki~\cite{kawasaki} had proved 
a general Riemann--Roch formula for orbifolds. His formula involves a 
sum over loci in $X$ with constant inertia group, and specifies to the above 
form in our case. Our proof relies on equivariant Riemann--Roch and a computation of 
intersection numbers on a resolution, and is independent 
of Kawasaki's; it gives the contributions in an explicit form (which could also 
be deduced from Kawasaki's result with extra work). Kawasaki's analytic formula was 
extended by Toen~\cite{toen} to the context of Deligne--Mumford 
stacks using algebraic methods.

Under further conditions on $(X,D)$, vanishing implies that 
$\chi(X,\bigo_X(nD))$ is simply the dimension of $H^0(X,\bigo_X(nD))$ for $n>0$.  
The finite dimensional vector spaces 
$H^0(X,\bigo_X(nD))$ fit together into a graded ring 
$$R(X,D)=\bigoplus_{n \geq 0} H^0(X,\bigo_X(nD)).$$ 
If $D$ is assumed ample as well as being $\bq$-Cartier, this ring is finitely 
generated. A surjection 
\[k[x_0, \ldots, x_n] \twoheadrightarrow R(X,D) \]
from a graded ring $k[x_0, \ldots, x_n]$ generated by variables $x_i$ of 
weights $a_i$ corresponds to an 
embedding $$i\colon X\cong\proj R(X,D)\hookrightarrow \bp(a_0,\ldots, a_n)$$ of~$X$ 
into a weighted projective space, with $\bigo_X(D)$ isomorphic 
to $\bigo_X(1)=i^{\ast}\bigo_{\bp}(1)$. 

Our aim is to construct quasi-smooth and well formed threefolds in weighted 
projective space. This implies in particular that the only singularities of $X$ are 
quotient singularities induced by the weights of the weighted projective space. 
Moreover, since we are interested in constructing threefolds with at worst 
canonical singularities, we restrict our study to curves of singularities which 
are generically of compound Du Val (cDV) type~\cite[Definition 2.1]{3-folds}. 

The main result of this paper is Theorem~\ref{thmdis}, which presents
an explicit RR formula of the above shape. As an application, we find new 
projective families of Calabi--Yau threefolds, via a study of their Hilbert 
series, which is brought into a compact form in Corollary~\ref{cyrhil}. 
The detailed and exhaustive analysis of families arising in this way, as well
as applications to other families of varieties such as Fanos and regular 
varieties of general type, will be presented elsewhere. 

\section{Definitions and notation}
\label{defnot}

We work over an algebraically closed field $k$ of characteristic zero. 
A $\bq$--\textit{divisor} on a normal variety $X$ is a formal linear combination 
of prime divisors with rational coefficients. 
A $\bq$--divisor $D$ is $\bq$--\textit{Cartier} if $mD$ is 
Cartier for some positive integer $m\neq 0$. In this case, 
if $C\subset X$ is a complete curve,
the \textit{degree} of $D$ on $C$ is defined by 
$$\deg D|_C=\frac{1}{m}\cdot \deg_C \bigo_X(mD).$$
If $X$ is projective, we also define the intersection $D\cdot c_2(X)$ as
$\frac{1}{m}f^{\ast}(mD)\cdot c_2(Y)$ computed on a resolution 
$f\colon Y\rightarrow X$, minimal over the cDV locus. 
A \textit{\cy threefold} is a normal projective threefold $X$ with canonical 
Gorenstein singularities, satisfying $K_X\sim \bigo_X$ (linear equivalence) 
and $H^1(X,\bigo_X)=0.$

A \textit{cyclic quotient singularity of type } $\frac{1}{r}(b_1,\ldots,b_n)$ is 
the quotient $\pi\colon \ba^n\rightarrow \ba^n/\mu_r$, where 
$\mu_r$ acts on $\ba^n$ by 
$$\mu_r\ni\varepsilon \colon\ (x_1,\ldots,x_n)\mapsto (\varepsilon^{b_1}x_1,\ldots,\varepsilon^{b_n}x_n)$$
We always assume that no factor of $r$ divides all the $b_i$, 
which is equivalent to the $\mu_r$-action being effective. 
The sheaf $\pi_{\ast}\bigo_{\ba^n}$ decomposes into eigensheaves
$$\bl_i=\{f\ | \ \varepsilon(f)=\varepsilon^i\cdot f\mbox{ for all } \varepsilon\in\mu_r\},$$
for $i=0,\ldots,r-1$. A singularity $Q\in X$ polarised by a Weil divisor $D$ is a
\textit{cyclic quotient singularity of type} 
$ _i\bigl(\frac{1}{r}(b_1,\ldots,b_n)\bigr),$
if $Q\in X$ is locally isomorphic to a point of type $\frac{1}{r}(b_1,\ldots,b_n)$ 
and $\bigo_X(D)\cong\bl_i$ near the singular point.

Let $(X,D)$ be a threefold $X$ containing a curve $C$ of singularities, 
equipped with a $\bq$-Cartier divisor $D$ which is Cartier away from $C$. 
Take a generic surface $S$ which intersects $C$ transversely 
in a finite number of points. Assume that every point in the 
intersection is a singular point of type $\ _k\bigl(\frac{1}{r}(1,-1)\bigr)$ 
on the polarised surface $(S, D|_S)$. Then $C\in X$ will be called 
a \textit{curve of singularities of transverse type} 
$\ _k\bigl(\frac{1}{r}(1,-1)\bigr)$, shorter a 
$\ _k\bigl(\frac{1}{r}(1,-1)\bigr)$ \textit{curve}, 
or often simply an $A_{r-1}$ \textit{curve}.

Note that every $A_{r-1}$ curve can contain a finite number 
of points of different type, 
which will be called \textit{dissident points.} A curve with dissident points is 
a \textit{dissident curve.} If an $A_{r-1}$ curve $C$ contains the 
dissident points $\{P_{\lambda}\ :\ \lambda\in\Lambda\}$ of types 
$\frac{1}{r \tau_{\lambda} }(a_{1\lambda},a_{2\lambda},a_{3\lambda})$, then 
we define the \textit{index of} $C$ to be 
\[\tau_C=\lcm_{\lambda\in\Lambda} \{\tau_{\lambda}\}.\] 
For an $A_{r-1}$ curve $C$ with no dissident points, $\tau_C=1.$

\begin{example}\rm The degree 13 hypersurface $X_{13}\subset \bp(1,1,2,3,6)$ in 
the weighted projective space with variables $x_1,x_2,y,z,t$ of the given degrees 
is a \cy threefold with two curves of singularities intersecting in a dissident 
point. The first $C_1=\{x_1=x_2=y=0\}$ is of type $\frac{1}{3}(1,2)$ and index~2 
because of the dissident point $(0,0,0,0,1)$ of type $\frac{1}{6}(1,2,3)$; 
likewise, $C_2=\{x_1=x_2=z=0\}$ is of type $\frac{1}{2}(1,1)$ and index~3.
\end{example}

For details on subvarieties in weighted projective space, consult~\cite{fletcher}. 
In particular, recall that a variety $X\subset\bp^{n}(a_0,\ldots,a_n)=\bp$ is 
\textit{quasi-smooth} if the affine cone $C_X\subset\ba^{n+1}$ over $X$ 
is smooth outside its vertex. In this case $X$ only has quotient singularities 
induced by the singularities of $\bp$. The weighted projective space 
$\bp(a_0,\ldots,a_n)$ is \textit{well formed} if 
$$(a_0,\ldots,\widehat{a_i},\ldots,a_n)=1 \ \mbox{ for each }i.$$ 
Moreover, $X$ of codimension $c$ in $\bp$ is \textit{well formed}, if it does 
not contain any $c+1$-codimensional singular stratum of $\bp$. 
Finally a pair $(X,D)$ of a variety $X$ and a $\bq$-divisor $D$ is quasi-smooth, 
respectively well formed, if there is an ample Cartier divisor $H$ on $X$, 
so that under the embedding $X\hookrightarrow\bp$ into a weighted projective 
space induced by the ample $\bq$-Cartier divisor $D+H$, the pair $X\subset\bp$ 
is quasi-smooth, respectively well formed.

\section{The Riemann--Roch formula}
\label{statethm} 

\subsection{The statement} 

\noindent Here is the main theorem of this paper:

\begin{theorem}
\label{thmdis}
Let~$(X,D)$ be a pair consisting of a normal projective threefold and a $\bq$-Cartier divisor, 
which is quasi-smooth and well formed. Assume further that the singularities 
of~$(X, D)$ consist of the following loci: 
\begin{itemize}
\item points $P\in X$ of type $_n\bigl(\frac{1}{s}(a_1,a_2,a_3)\bigr)$ (dissident and 
isolated), and 
\item curves $C\subset X$ of generic type $_k\bigl(\frac{1}{r}(1,-1)\bigr)$ with index 
$\tau_C$. 
\end{itemize}
Then for all positive integers $m$,
$$\begin{array}{rcl}
\chi(X,\bigo_X(mD)) &= & \chi(\bigo_X)+\frac{1}{12}mD(mD-K_X)(2mD-K_X)+m\displaystyle\frac{D\cdot c_2(X)}{12} \\
 & & \\
  & & \mbox{}+\displaystyle\sum_{P}c_P(mD)+\sum_{C}s_C(mD),
\end{array}$$ 
where
$$c_P(mD)=
\displaystyle\frac{1}{s}\sum_{\substack{\ep\in\mu_s \\ \ep^{a_i}\ne1\,\forall i = 1,2,3}}
\frac{\varepsilon^{-nm}-1}{(1-\varepsilon^{-a_1})(1-\varepsilon^{-a_2})(1-\varepsilon^{-a_3})}$$
and
$$\begin{array}{rcl}
s_C(mD)& = & -m\displaystyle\frac{\ol{mk}(r-\ol{mk})}{2r}\deg D|_C + \frac{\ol{mk}(r-\ol{mk})}{4r}\deg K_X|_C \\
 & & \\
      & & +\displaystyle\frac{\ol{mk}(r-\ol{mk})(r-2\cdot\ol{mk})}{12r^2\tau_C}N_C,   
\end{array}$$
where the integer $N_C$ is an invariant of $X$ in a neighbourhood of~$C$, and 
$\overline{\phantom{\Sigma}}$ denotes the smallest residue mod~$r$. 
\end{theorem}

\begin{remark}\rm
Note that a point $P\in X$ of type $_n\bigl(\frac{1}{s}(a_1,a_2,a_3)\bigr)$ is either
\begin{enumerate}
\item an isolated point singularity if $\hcf(a_i,s)=1$ for all $i=1,2,3,$\ \ or 
\item a dissident point on some curve if $\hcf(a_i,s)=\alpha_i\neq 1$ for some $i$. 
\end{enumerate}
In the first case, $c_P(mD)$ equals to the so-called basket contribution 
to Riemann--Roch~\cite{ypg}. In the second case, as the singularities along curves are 
of $A_{\alpha_i-1}$ transverse type, we must have $a_j+a_k=0$ mod $\alpha_i$ for different 
indexes $i,j,k$. Since $X$ is well formed, $a_i,a_j,s$ have no common divisor.
\end{remark}

\subsection{The outline of the proof}
\label{outline}

Let $(X,D)$ be a pair satisfying the conditions of Theorem~\ref{thmdis}. 
Choose a projective resolution $f\colon Y\rightarrow X$ which is crepant at the generic point of 
each one-dimensional component of the singular locus of $(X,D)$. The sheaf $\bigo_X(D)$ is 
a rank-1 reflexive sheaf on $X$; define $\L= f^*(\bigo_X(D))/({\rm torsion})$, 
a rank-one torsion-free sheaf on $Y$. Let also $H$ be an ample Cartier divisor on $X$. 

For some integers $n,N$, there exists a surjection 
$\bigo_X^N(-nH)\twoheadrightarrow\bigo_X(D)$ on $X$ 
which pulls back to a surjection 
$\bigo_Y^N\otimes f^*\bigo_X(-nH)\twoheadrightarrow\L$ on $Y$ 
and hence gives an exact sequence
$$0\to \bk \to \bigo_Y^N\otimes f^*\bigo_X(-nH) \to\L \to 0.$$
Under $f_\ast$, this becomes the long exact sequence
$$\begin{array}{cccccccc}
0\to & f_{\ast}\bk    & \to & \bigo_X^N(-nH) & {\to} & f_{\ast}\L  & \to \\[6pt]
                 & R^1f_{\ast}\bk & \to & 0 &
\to                   & R^1f_{\ast}\L & \to \\[6pt]
                 & R^2f_{\ast}\bk &\to  & 0 &
\to                   & R^2f_{\ast}\L & \to  & 0,
\end{array}$$
where we used $f_{\ast}\bigo_Y\cong\bigo_X$ which holds since $X$ is normal, and 
the projection formula together with $R^if_{\ast}\bigo_Y=0$ for $i>0$ which holds 
as quotient singularities in characteristic zero are rational. It is easy to
see that $f_\ast\L\cong\bigo_X(D)$; from the exact sequence we also have 
$R^2f_{\ast}\L=0$ and $R^1f_{\ast}\L\cong R^2f_{\ast}\bk$, the latter 
necessarily supported on the isolated and dissident singularities of $(X,D)$.

The reflexivization of the rank-one torsion-free sheaf $\L$ on the smooth variety
$Y$ is a line bundle $\bigo_Y(D_Y)$. As $\bigo_X(D)$ is saturated, we also have
$f_\ast\bigo_Y(D_Y)\cong\bigo_X(D)$. Thus the exact sequence
\[ 0 \to \L \to \bigo_Y(D_Y) \to \Q \to 0\]
give rise to an injection
\[ f_\ast \Q \hookrightarrow R^1f_\ast \L \cong R^2 f_\ast \bk, \] 
proving that $f_\ast Q$ is supported on the isolated and dissident singular points 
of $(X,D)$. As $\Q$ is supported in codimension two on $Y$, this implies that all its 
higher pushforwards are also supported on these singular points.

Finally, by the Leray spectral sequence
\[
\chi(Y, \L) = \sum_i (-1)^i \chi(X, R^if_{\ast}\L), 
\]
which can be rewritten, using the above exact sequences and isomorphisms, as
\begin{eqnarray}\label{firststep}
\chi(X,\bigo_X(D)) & = & \chi(Y, \bigo_Y(D_Y)) + \chi(X, R^1f_{\ast}\bl) - \sum_i (-1)^i \chi(X, R^if_\ast \Q)\nonumber\\
& = & \chi(Y, \bigo_Y(D_Y)) + \P_1,
\end{eqnarray}
where, by our earlier remarks, $\P_1$ is a contribution from sheaves supported entirely 
on the isolated and dissident singular points of $(X,D)$. Our arguments in fact imply 
that this contribution is local in a stronger sense: it only depends on the analytic type 
of the isolated and dissident quotient singularities of $(X,D)$. 
This holds since the constructions of $\L$ and $\Q$ are universal, and an analytic 
isomorphism preserving the type of the quotient singularities 
necessarily gives an analytic isomorphism between these sheaves, and 
thus an equality of Euler characteristics. 

In the next step we will express $D_Y$ and $K_Y$ in terms of $D, K_X$ and 
the exceptional divisors of the resolution. For a singular point $P\in X$, let $\{{}^P\! G_j\}$
be the exceptional surfaces mapping to $P$ under $f$; similarly, for a curve $C\subset X$ 
of singularities, let $\{{}^C\! E_i\}$ and be the exceptional surfaces
mapping surjectively to $C$. Note that every dissident or isolated singular point $P$ 
is locally analytically isomorphic to $\ba^3/\mu_s$, and the configuration 
of $\{{}^P\! G_j\}$ depends only on the analytic singularity type of $P$. 

As $f$ is crepant at the generic point of each curve $C$, we have 
$$K_Y=f^{\ast}K_X+N,\mbox{ where } N=\sum_P\sum_j\,^P\!\gamma_{j}\, ^P\!G_j;$$
also
$$f^{\ast}D=D_Y+\sum_{C}R_C+M, \mbox{ where } R_C=\sum_i\,^C\!\alpha_{i}\,^C\!E_i 
\mbox{ and } M=\sum_P\sum_j\,^P\! \beta_{j}\,^P\!G_j,$$
with ${}^C\!\gamma_j,{}^C\!\alpha_i, {}^P\!\beta_j\in\bq$. Here $f^{\ast}D$ is by 
definition $\frac{1}{m}f^{\ast}(mD)$ for an integer $m$ which makes $mD$ Cartier.  

\begin{lemma}
\label{lemconfig} If $\overline{\phantom{\Sigma}}$ denotes the smallest residue mod~$r$, then
$$R_C=\displaystyle\frac{k}{r}\,^C\!E_1+\frac{\overline{2k}}{r}\, ^C\!E_2+\ldots+\frac{\overline{(r-1)k}}{r}\, ^C\!E_{r-1}.$$
\end{lemma}

\begin{prf} Choose a general transverse hypersurface $S$ and a point $Q\in S\cap C$. 
By assumption, on the polarised $(S,D|_S)$ the type of 
$Q$ is $_k\bigl(\frac{1}{r}(1,-1)\bigr)$.
Such a point has a unique crepant resolution 
 $f|_{\tilde{S}}\colon \tilde{S}\rightarrow S$, and we necessarily have 
$$(f|_{\tilde{S}})^\ast D|_S=D_Y|_{\tilde{S}}+
\frac{k}{r}F_1+\frac{\overline{2k}}{r}F_2+\ldots+
\frac{\overline{(r-1)k}}{r}F_{r-1},$$
where the $F_i$ are exceptional lines in $\tilde S$ forming an~$A_{r-1}$ 
configuration. 

Each $F_i$ is a fibre in one of the surfaces $^C\!E_j$. Since $X$ is assumed quasi-smooth 
and well formed, there is no ramification in codimension 1; thus $k$ and $r$ are coprime.
Therefore the integers $k, \overline{2k},\ldots,\overline{(r-1)k}$ 
are all different, which implies that $F_1,\ldots,F_{r-1}$ must be fibres in different 
surfaces $^C\!E_1,\ldots,^C\!E_{r-1}$. Indeed, if some $F_i$ would be fibres in the same
surface, the relevant coefficients should be equal, arising from an irreducible
divisor on the threefold. The statement follows. 
\end{prf}
\ethrm

\begin{figure}[ht] 
\begin{center}
\input{conf.pstex_t}
\caption{Resolution of an $A_{r-1}$ curve $C$ with a dissident point $P$}
\label{pict5}
\end{center}
\end{figure}
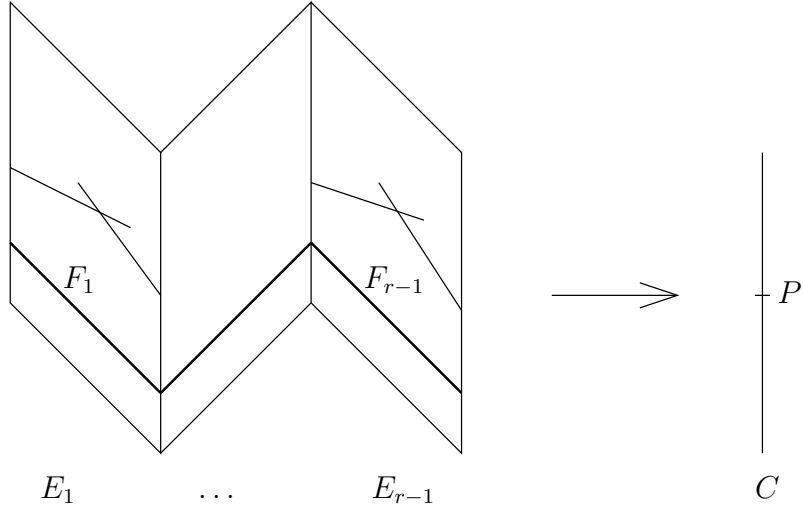

In particular, there are $(r-1)$ irreducible exceptionals
${^C\!E_1},\ldots,{^C\!E_{r-1}}$ over every $C$, each isomorphic to a blowup in a few points 
(over the dissident points) of a smooth surface geometrocally ruled over $C$.   
Moreover,  
\label{conref}
$^C\!E_i$ and $^C\!E_j$ intersect along a curve isomorphic to $C$ if $j=i+1$ and are 
disjoint otherwise. The situation is shown on Figure~\ref{pict5}. For two different curves $C$ 
and $\bar{C}$, the surfaces $^C\!E_i$ and $ ^{\bar{C}}\!E_j$  intersect only over the 
dissident points.

We proceed to show that many intersection numbers between these divisors actually vanish. 
Choose an ample Cartier divisor $H$ on $X$ such that $D+H$ is also ample. 
Take $t$ large and divisible, so that the linear systems $|\,t\,(D+H)|$ and $|tH|$
contain nonsingular divisors $S_1$ and $S_2$ respectively.
Then $f^{\ast}S_1\cdot{^C\!E_i}$ is a finite union of fibres in $^C\!E_i$, since
$f^{\ast}S_1|_{^C\!E_i}$ is a pullback of $t(D+H)|_{C}$; the same holds for
$f^{\ast}S_2\cdot{^C\!E_i}$. Then 
\begin{eqnarray*} f^{\ast}(tD)\cdot f^{\ast}(tD)\cdot {^C\!E_i}&=&(f^{\ast}S_1-f^{\ast}S_2)\cdot (f^{\ast}S_1-f^{\ast}S_2) \cdot {^C\!E_i}\\ 
& = & f^{\ast}S_1\cdot (\mbox{fibres}) - f^{\ast}S_2\cdot (\mbox{fibres})\\ 
& = & S_1\cdot f_{\ast}(\mbox{fibres}) - S_2\cdot f_{\ast}(\mbox{fibres})\\
&=&0.
\end{eqnarray*}
using also the projection formula. Thus
$f^{\ast}D\cdot f^{\ast}D\cdot R_C=0$
and similarly $f^{\ast}D\cdot f^{\ast}K_X\cdot R_C=0.$
The projection formula also gives
$$(f^{\ast}D)^3=D^3,\ (f^{\ast}D)^2\cdot f^{\ast}K_X=D^2K_X,\ f^{\ast}D\cdot(f^{\ast}K_X)^2=DK_X^2.$$
The intersection numbers involving $M$ or $N$
depend only on the isolated or dissident singular points of $X$.
Moreover, for different curves of singularities, the intersections 
$$f^{\ast}D\cdot R_C R_{\bar{C}},\ f^{\ast}K_X\cdot R_C R_{\bar{C}}\ \mbox{ and }\ (R_C)^2R_{\bar{C}}$$
can only be nonzero because of intersections over dissident points.

Finally using~(\ref{firststep}), Riemann--Roch for the smooth threefold $Y$, and the birational
invariance of $\chi(\bigo_X)$, we obtain
\begin{eqnarray}
\label{eqena}
\chi(X, \bigo_X(D))& = & \chi(\bigo_X)+\frac{1}{12}D(D-K_X)(2D-K_X)+\frac{1}{12}f^{\ast}D\cdot c_2(Y) + \nonumber\\
& &  \frac{1}{2}f^{\ast}D \cdot \sum_CR_C^2-\frac{1}{4}f^{\ast}K_X \cdot \sum_CR_C^2 -\nonumber\\
& &  \frac{1}{6}\sum_CR_C^3-  \frac{1}{12}\sum_CR_C\cdot c_2(Y)+ \P_2. 
\end{eqnarray}
Here $\P_2$ denotes a quantity that depends, in addition to the earlier quantity $\P_1$, on 
intersection numbers of divisors over the isolated and dissident singularities of $(X,D)$. 
Since the latter are purely analytical, $\P_2$ also depends only on the analitic type
of these singularities. 

\subsection{The contribution from a curve of singularities} 

In this section we focus on the part of the contribution that depends only on a small 
analytic neighbourhood of the curves of singularities. During the argument we will
often meet divisors and intersection numbers which depend on the configuration of 
exceptional divisors over the dissident points, and thus depend only on the analytic 
type of these points. We will denote all such divisors by $\D$ and numerical 
contributions by $\P$, and will not worry about the exact expressions. 

The next two lemmas will analyze the intersection numbers in~(\ref{eqena}).

\begin{lemma}
\label{lemena}
For a $\ _k\bigl(\frac{1}{r}(1,-1)\bigr)$ curve $C$ of singularities of $(X,D)$,  
$$\frac{1}{2}f^{\ast}D\cdot R_{C}^2=-\deg D|_C\cdot \frac{\ol{k}(r-\ol{k})}{2r}$$
and  
$$\frac{1}{4}f^{\ast}K_X\cdot R_{C}^2=-\deg K_X|_C\cdot \frac{\overline{k}(r-\overline{k})}{4r}. $$
\end{lemma}

\begin{prf}
If $S\subset X$ is surface which intersects $C$ transversely in all points, then
the points in $C\cap S$ are all of type $\ _k\bigl(\frac{1}{r}(1,-1)\bigr)$ on $(S,D|_S)$. 
Riemann--Roch for a resolution $f|_{\tilde{S}}\colon \tilde{S}\rightarrow S$ of these 
points gives
$$  \chi(S, \bigo_S(D))=\chi(\bigo_S)+\frac{1}{2}D|_S(D|_S-K_S)+\frac{1}{2}B_{\ol{k}}^2,$$
as in the proof of Theorem 9.1 in \cite{ypg}. Here $B_{\ol{k}}=R_C|_{\tilde{S}}$
is the $\bq$-divisor on the surface $\tilde S$ contributing the so-called basket contribution 
$$\frac{1}{2}B_{\ol{k}}^2 = \#(C\cap S)\cdot\left(-\frac{\overline{k}(r-\overline{k})}{2r}\right).$$
As before, for some integer $t$ we can write $tD$ as a difference of two 
smooth surfaces $S_1, S_2$ both intersecting $C$ transversally; thus
$$
\frac{1}{2}f^{\ast}D\cdot R_{C}^2=\frac{1}{2t}(f^{\ast}S_2-f^\ast S_1)\cdot R_{C}^2=  
-\deg D|_C\cdot \frac{\overline{k}(r-\overline{k})}{2r}.
$$
A same argument also shows
$$\frac{1}{4}f^{\ast}K_X\cdot R_{C}^2=-\deg K_X|_C\cdot \frac{\overline{k}(r-\overline{k})}{4r}$$
which ends the proof of Lemma \ref{lemena}.
\end{prf}
\ethrm

\begin{lemma}
\label{lemdva}
Let $C$ be a $\, _k\bigl(\frac{1}{r}(1,-1)\bigr)$ curve in $(X,D)$. Then
$$\begin{array}{l}
R_{C}^3 + \frac{1}{2}R_{C}\cdot c_2(Y)=  \\
\\
-\frac{1}{r^2}\ol{k}(r-\ol{k})(r-2\ol{k})
 \left[(r-2)\left(1-g+\frac{1}{2}K_Y\,^C\!E_1\,^C\!E_2\right)+\,^C\!E_{r-2}^2\,^C\!E_{r-1}\right]+\P,
\end{array}$$
where $g$ is the genus of $C$.
\end{lemma}

\begin{prf}
>From the definition of $R_C$,
$$
R_{C}^3 =  \sum_{i=1}^{r-1} \left(\displaystyle\frac{\ol{ik}}{r}\right)^3E_i^3+3\sum_{i=1}^{r-2} \left(\displaystyle\frac{\ol{ik}}{r}\right)^2\frac{\ol{(i+1)k}}{r}E_i^2E_{i+1}+
\frac{\ol{ik}}{r}\left(\frac{\ol{(i+1)k}}{r}\right)^2E_iE_{i+1}^2
$$
and
$$R_{C}\cdot c_2(Y)  =\left(\frac{\ol{k}}{r}E_1+\frac{\ol{2k}}{r}E_2+\ldots+\frac{\ol{(r-1)k}}{r}E_{r-1}
\right)\cdot c_2(Y).$$
We will simplify these expressions using the properties of ruled surfaces.
Recall that each $E_i$ is isomorphic to a smooth ruled surface $\hat{E}_i$ blown up in a number 
of points. The blowups $E_i\stackrel{\pi_i}{\to}\hat{E}_i$ happen only over the 
dissident points. In other words, the resolution
$f\colon Y\rightarrow X$ contracts Exc$(\pi_i)\subset E_i$ to the dissident points on~$C$.

Denote $\gamma_i=E_i\cap E_{i+1}$ which is isomorphic to $C$ of genus $g.$
Note that $\pi_i(\gamma_i)=\pi_{i+1}(\gamma_i)$ is a section in both ruled surfaces $\hat{E}_i$ and $\hat{E}_{i+1}$. 
Then
$$K_{E_i}=\pi_i^{\ast}K_{\hat{E}_i}+ \D \sim -2\gamma_i+(2g-2+\gamma_i^2)f_i+\D,$$ 
where $f_i$ is a generic fibre on $E_i$. In particular 
$K_{E_i}^2=8(1-g)+\P.$ Similarly
$$K_{E_{i+1}}\sim -2\gamma_i+(2g-2+\gamma_i^2)f_{i+1}+\D.$$ 
By the adjunction formula we can compute
\begin{eqnarray*}
(\gamma_i^2)_{E_{i+1}}& = & (K_{E_i}-K_Y|_{E_i})E_{i+1}|_{E_i}\\
& = & -(\gamma_i^2)_{E_i}+2(g-1)-K_YE_iE_{i+1}+\P, 
\end{eqnarray*}
where $(\gamma^2)_E$ is the self-intersection of the curve $\gamma$ computed in the surface
$E\supset\gamma$. 
Similarly 
$$(\gamma_i^2)_{E_i}=-(\gamma_i^2)_{E_{i+1}}+2(g-1)-K_YE_iE_{i+1}+\P.$$
Together these give \label{rulled}
$$E_i^2E_{i+1}+E_iE_{i+1}^2=2(g-1)-K_YE_iE_{i+1}+\P\ \mbox{ for }i=1,\ldots,r-2.$$

Next, 
$$E_i^2E_{i+1}+E_{i+1}E_{i+2}^2=(\gamma_i^2)_{E_{i+1}}+(\gamma_{i+1}^2)_{E_{i+1}}+\P
\ \hbox{ for }\ i=1,\ldots,r-3$$
since $\pi_{i+1}(\gamma_i)$ and $\pi_{i+1}(\gamma_{i+1})$ are disjoint 
sections in~$\hat{E}_{i+1}.$

By the projection formula also the following holds:
$$K_YE_{i-1}E_{i}-K_YE_iE_{i+1}=\P\mbox{ and }K_YE_i^2+2K_YE_iE_{i+1}=\P.$$

Finally, putting together all of the above yields
$$\begin{array}{rcl}
E_i^2E_{i+1} & =&-(r-2-i)2(g-1)+(r-2-i)K_YE_1E_2+E_{r-2}^2E_{r-1}+\P,  \\[6pt]
E_iE_{i+1}^2 & =&(r-1-i)2(g-1)-(r-1-i)K_YE_1E_2-E_{r-2}^2E_{r-1}+\P, 
\end{array}$$
for all $i=1,\ldots,r-2$. 

Expressions involving the second Chern class $c_2(Y)$ can also be computed without 
difficulty. Let $\bn_i$ be the normal bundle of $E_i$ in 
the smooth threefold $Y$, and $\bt_{E_i},\ \bt_Y$ the tangent sheaves. 
There is an exact sequence
$$0\to \bt_{E_i}\to j_i^{\ast}\bt_Y \to \bn_i \to 0,$$
where $j_i\colon E_i\rightarrow Y$ is the inclusion.
Comparison of Chern polynomials and adjunction gives
$$
c_2(Y)\cdot E_i= c_2(E_i)-K_{E_i}E_i|_{E_i} = c_2(E_i)-(K_Y+E_i)E_i^2. 
$$ 

Also $c_2(E_i)$ can be expressed with other invariants
$$\frac{1}{12}\left(K_{E_i}^2+c_2(E_i)\right)=1+p_a(E_i)=1-g.$$
This in particular implies $c_2(E_i)=4(1-g)+\P.$
From
$$8(1-g)+\P=K_{E_i}^2=(K_Y+E_i)^2E_i=E_i^3+2K_YE_i^2+\P$$
we get
$$c_2(Y)\cdot E_i=-\frac{1}{2}E_i^3+\P=4(g-1)-2K_YE_1E_2+\P.$$

This reduces the expression for $R_C^3$ to
$$\begin{array}{l}
\displaystyle\frac{8(1-g)+4K_YE_1E_2}{r^3}\sum_{i=1}^{r-1}\left(\ol{ki}\right)^3+ \\
\displaystyle\frac{3}{r^3}\sum_{i=1}^{r-2} (\ol{ik})^2\ol{(i+1)k}\Big(-2(r-2-i)(g-1)+ (r-2-i)K_YE_1E_2+ E_{r-2}^2E_{r-1}\Big)+\\
\displaystyle\frac{3}{r^3}\sum_{i=1}^{r-2} \ol{ik}(\ol{(i+1)k})^2\Big(2(r-1-i)(g-1)- (r-1-i)K_YE_1E_2- E_{r-2}^2E_{r-1}\Big)+ \P 
\end{array}$$
and $R_{C}\cdot c_2(Y)$ reduces to
$$-\frac{4(1-g)+2K_YE_1E_2}{r}\sum_{i=1}^{r-1} \ol{ki}+\P.$$

Arguing as in \cite{morrison}, we can simplify 
$R_{C}^3 + \frac{1}{2}R_{C}\cdot c_2(Y)$ to the form
$$-\frac{1}{r^2}\ol{k}(r-\ol{k})(r-2\ol{k})
\left[(r-2)\left(1-g+\frac{1}{2}K_YE_1E_2\right)+E_{r-2}^2E_{r-1}\right]+\P.$$
This ends the proof of Lemma~\ref{lemdva}.
\end{prf}
\ethrm

\begin{corollary}
\label{rrformm}
Let $(X,D)$ satisfy the conditions of Theorem~\ref{thmdis}. Then for all positive integers $m$, 
\begin{eqnarray*}
\chi(X,\bigo_X(mD))&= & \chi(\bigo_X)+\frac{1}{12}mD(mD-K_X)(2mD-K_X) + m\frac{D\cdot c_2(X)}{12} +\\
 & &\\
 & &\sum_{C}s_C(mD)+ \P_3
\end{eqnarray*} 
where
$$\begin{array}{c}
s_C(mD) =-m\displaystyle\frac{\ol{mk}(r-\ol{mk})}{2r}\deg D|_C + \frac{\ol{mk}(r-\ol{mk})}{4r}\deg K_X|_C + \\
 \\
  \displaystyle\frac{\ol{mk}(r-\ol{mk})(r-2\cdot \ol{mk})}{6r^2}\left( (r-2)(1-g)+
 \frac{r-2}{2}K_Y\,^C\!E_1\,^C\!E_2+\,^C\!E_{r-2}^2\,^C\!E_{r-1}\right),   
\end{array}$$
and $\P_3$ is a contribution from the dissident and isolated singular points of $(X,D)$. 
\end{corollary}

\begin{prf}
For $m=1$, we only need to put the results of Lemma~\ref{lemena} and
Lemma~\ref{lemdva} into Formula~(\ref{eqena}) on page~\pageref{eqena}.
For $m>1$, consider $\bigo_X(mD)$, which is of transverse type 
$_{\ol{mk}}\bigl(\frac{1}{r}(1,-1)\bigr)$ on $C$, and repeat the proof using
$$f^{\ast}(mD)=D^{(m)}+\sum_{C}R_C^{(m)}+\D,$$
where $D^{(m)}$ is a Cartier divisor on $Y$, and 
$$R_C^{(m)}=\frac{\overline{mk}}{r}\,^C\!E_1+\frac{\overline{2mk}}{r}\, ^C\!E_2+\ldots+
\frac{\overline{(r-1)mk}}{r}\,^C\!E_{r-1}.$$
\end{prf}
\ethrm

\subsection{The contribution from a dissident or isolated singular point}
\label{three}

The analytically invariant contributions to RR from dissident 
or isolated singular points, denoted by $\P_3$ in Corollary~\ref{rrformm},
can be computed on any model that contains such singularities. 
We begin by showing the existence of such projective varieties.

\begin{proposition}
\label{exag}
Fix positive integers $s$ and $a_1,a_2,a_3$, and assume that $a_i,a_j,s$ have no common factor 
for all different $i,j\in\{1,2,3\}$. 
There exists a smooth projective 
3--fold $Z$ together with an action of $\mu_s$ with the following properties: 
the action fixes a number of points on which a generator $\varepsilon\in\mu_s$ acts by 
$$\varepsilon\colon\ z_1,z_2,z_3\ \mapsto\ \varepsilon^{a_1}z_1,\varepsilon^{a_2}z_2,\varepsilon^{a_3}z_3.$$
If $\hcf(a_i,s)=\alpha_i\neq 1$, these points lie on curves
which are fixed by $\varepsilon^{\frac{s}{\alpha_i}}.$
Finally $\varepsilon^{\frac{s}{\alpha_i}}\in\mu_{\alpha_i}$ acts in the normal direction 
of each curve by
$$\varepsilon^{\frac{s}{\alpha_i}}\colon\ z_j,z_k\ \mapsto\ \varepsilon^{\frac{s}{\alpha_i}a_j} z_j,
\varepsilon^{\frac{s}{\alpha_i}a_k}z_k,\ \mbox{ for } j,k\in\{1,2,3\}-\{i\}$$
and freely away from the curve. 
\end{proposition}

\begin{prf}
We imitate the proof of (8.4) in~\cite{ypg}. Choose an integer $l\geq 3$ and consider action of
$\mu_s$ on $\bp^{l+3}(1,1,\ldots,1)$ given by 
$$x_1,x_{2},x_{3},x_{4},\ldots,x_{l+4}\ \mapsto\ 
\varepsilon^{a_1} x_{1},\varepsilon^{a_2} x_{2},\varepsilon^{a_3}x_{3},x_{4},\ldots,x_{l+4}.$$
This action fixes $\bp^{l}=\{x_1=x_2=x_3=0\}$ and acts in the normal direction by
$$\ x_1,x_2,x_3\ \mapsto\ \varepsilon^{a_1}x_1,\varepsilon^{a_2}x_2,\varepsilon^{a_3}x_3.$$
If $\hcf(a_i,s)=\alpha_i\neq 1$, the action is not free on $\bp^{l+1}=\{x_j=x_k=0\}$ for $j,k\in\{1,2,3\}-\{i\}$. 
This is fixed by 
$\varepsilon^{\frac{s}{\alpha_i}}$ which acts in the normal direction by
$$x_j,x_k\ \mapsto\ \varepsilon^{\frac{s}{\alpha_i}a_j} x_j,
\varepsilon^{\frac{s}{\alpha_i}a_k}x_k.$$
Another locus on which the action 
might not be free, is $\{x_{4},\ldots,x_{l+4}=0\}$. We will avoid this locus by defining
$$X\subset\bp^{l+3}/\mu_s,$$
as a complete intersection of $l$ general very ample divisors. Let $Z$ be the inverse image of 
$X$ under the quotient $\bp^{l+3}\rightarrow \bp^{l+3}/\mu_s $. 
Such $Z$ clearly satisfies the conditions in the proposition. 
\end{prf}
\ethrm

Let $X$ be a projective threefold with a singularity of type $\frac{1}{s}(a_1,a_2,a_3)$ 
as described in Proposition~\ref{exag}. 
We restricted our study to Du Val curves, which means that if
$\hcf(a_i,s)=\alpha_i\neq 1,$ then
$$a_j+a_k=0\ \mbox{ mod }\ \alpha_i.$$ 
Thus $\sing X$ equals to a number of dissident 
$$A_{\alpha_i-1}=\frac{1}{\alpha_i}(a_j,a_k)$$ curves, which intersect in dissident 
points of type
$\frac{1}{s}(a_1,a_2,a_3).$

Write $\pi\colon Z\rightarrow X$ for the quotient map 
and let $\bl_{n}$ be the $n$th eigensheaf of the action of 
$\varepsilon\in\mu_s$ on $\pi_{\ast}\bigo_Z.$ Then
$$\pi_{\ast}\bigo_Z=\bigoplus_{n=0}^{s-1}\bl_{n}$$
implies
$$H^p(Z,\bigo_Z)=\bigoplus_{n=0}^{s-1}H^p(X,\bl_{n}).$$
The group action on any $f\in\bl_{n}$ is given by $\varepsilon (f)=
\varepsilon^n \cdot f$, thus
$$\mbox{Tr}(\varepsilon : H^p(Z,\bigo_Z))=
\sum_{n=0}^{s-1}h^p(X,\bl_{n})\cdot \varepsilon^{n}. $$
Moreover,
$$\sum_p(-1)^p\mbox{Tr}(\varepsilon : H^p(Z,\bigo_Z))=
\sum_{n=0}^{s-1}\chi(X,\bl_{n})\cdot \varepsilon^{n}. $$
In order to simplify the notation, denote $\sum_p(-1)^p\mbox{Tr}(\varepsilon : H^p(Z,\bigo_Z))$ by $A_{\varepsilon}.$
Then
$$\sum_{n=0}^{s-1}\chi(X,\bl_{n})\cdot \varepsilon^{n}=A_{\varepsilon} \ \ \mbox{ and } \ \ \sum_{n=0}^{s-1}\chi(X,\bl_{n})=\chi(\bigo_Z).$$

The last two formulas can be considered as a linear system of $s$ equations in $\chi(X,\bl_{n})$ and 
variable $\varepsilon$, 
$$\left(\begin{array}{ccccc}
1 &  1            &     1         & \ldots & 1 \\
1 &  \varepsilon  & \varepsilon^2 & \ldots & \varepsilon^{s-1}  \\
1 & \varepsilon^2 & \varepsilon^4 & \ldots &  \varepsilon^{2(s-1)}  \\
\ldots  & \ldots  &  \ldots       & \ldots &  \\
1 & \varepsilon^{s-1} & \ldots    &  \ldots  & \varepsilon 
\end{array}\right)
\left(\begin{array}{c}
\chi(X,\bl_{0})\\
\chi(X,\bl_{1})\\
\chi(X,\bl_{2})\\
\vdots \\
\chi(X,\bl_{s-1})
\end{array}\right)=
\left(\begin{array}{c}
\chi(\bigo_Z)\\
A_{\varepsilon}\\
A_{\varepsilon^2}\\
\vdots\\
A_{\varepsilon^{s-1}}\\
\end{array}\right).$$
Write $\chi(X,\bl_{0})=\chi(\bigo_X)$ and eliminate $\chi(\bigo_Z)$ from the 
solution. We end up with 
\begin{eqnarray}
\label{eqypg}
\chi(X,\bl_{n}) &=&\chi(\bigo_X)+\frac{1}{s}\sum_{j=1}^{s-1}(\varepsilon^{-jn}-1)A_{\varepsilon^j}  \\
&=&\chi(\bigo_X)+\frac{1}{s}\sum_{\varepsilon\in\mu_s}(\varepsilon^{-n}-1)A_{\varepsilon}. \nonumber
\end{eqnarray}

We can compute $A_{\varepsilon}$ for all $\varepsilon\in\mu_s$ using the 
Atiyah--Singer--Segal equivariant RR formula~\cite{as},~\cite[p.565]{atiyah}.
If $\varepsilon^{a_i}\neq 1$ for all $i=1,2,3$, then the fixed locus of 
$\varepsilon$ is a number of points. In this case, for each point
$$
A_{\varepsilon}=\sum_p(-1)^p\mbox{Trace}(\varepsilon|H^p(Z,\bigo_Z)) = 
\frac{1}{(1-\varepsilon^{-a_1})(1-\varepsilon^{-a_2})(1-\varepsilon^{-a_3})}.
$$
If, on the other hand $\varepsilon^{a_i}= 1$ for one of the $i=1,2,3,$
then the fixed locus of $\varepsilon$ is a curve $C_i$. In this case, 
$\hcf(a_i,s)=\alpha_i\neq 1$ and $\varepsilon$ is an element of $\mu_{\alpha_i}$. 
The action in the normal direction equals 
$$\varepsilon^{a_j},\varepsilon^{a_k}=\varepsilon^{a_j},\varepsilon^{-a_j}.$$
Thus, we can use the equivariant Riemann--Roch formula for the cyclic group 
$\mu_{\alpha_i}$. Let $C_i$ have genus $g_i$, and let $x_1,x_2$ denote the first 
Chern classes of $\bn_{C_i/Z}=\bn_1\oplus\bn_2$. We get
$$\begin{array}{rcl}
A_{\varepsilon} & = & \sum_p(-1)^p\mbox{Trace}(\varepsilon|H^p(Z,\bigo_Z)) \\
 & & \\
 & = & \deg\left\{\displaystyle\frac{\mbox{td}(\bt_C)}{(1-\varepsilon^{-a} e^{-x_1})(1-\varepsilon^{-(r-a)} e^{-x_2})}\right\}_1\\
 & & \\
 & = & \frac{1-g_i}{(1-\varepsilon^{-a_j})(1-\varepsilon^{a_j})}
-\frac{\varepsilon^{-a_j}}{(1-\varepsilon^{-a_j})^2(1-\varepsilon^{a_j})}\deg_{C_i} x_1 
-\frac{\varepsilon^{a_j}}{(1-\varepsilon^{-a_j})(1-\varepsilon^{a_j})^2}\deg_{C_i} x_2. 
\end{array}$$

Since we assumed that $a_i,a_j,s$ have no common divisor for all $i,j=1,2,3$, 
we can split the sum over $\{\varepsilon\in\mu_s\}$ into four 
subsumes over (possibly empty) disjoint sets 
$$\begin{array}{l}
\{\varepsilon\in\mu_s \mbox{ s.t. }\varepsilon^{a_i}\neq 1\ \forall i=1,2,3\}, \\
\\
\{\varepsilon\in\mu_s \mbox{ s.t. }\varepsilon^{a_1}=1\}=\{\varepsilon\in\mu_{\alpha_1}\},\\
\\
\{\varepsilon\in\mu_s \mbox{ s.t. }\varepsilon^{a_2}=1\}=\{\varepsilon\in\mu_{\alpha_2}\} \mbox{ and}\\
\\
\{\varepsilon\in\mu_s \mbox{ s.t. }\varepsilon^{a_3}=1\}=\{\varepsilon\in\mu_{\alpha_3}\}.
\end{array}$$
This rewrites Formula~(\ref{eqypg}) into
$$\begin{array}{l}
\chi(X,\bl_{n})= \chi(\bigo_X) +\displaystyle\frac{1}{s}\sum_{\substack{\ep\in\mu_s \\ \ep^{a_i}\ne1\,\forall i = 1,2,3}} 
\frac{\varepsilon^{-n}-1}{(1-\varepsilon^{-a_1})(1-\varepsilon^{-a_2})(1-\varepsilon^{-a_3})} + \\
\\
\displaystyle\frac{1}{s}\sum_{(i,j)}\,\sum_{\scriptscriptstyle\varepsilon\in\mu_{\alpha_i}}(\varepsilon^{-n}-1)
\left( \frac{1-g_i}{(1-\varepsilon^{-a_j})(1-\varepsilon^{a_j})}
-\frac{\varepsilon^{-a_j}}{(1-\varepsilon^{-a_j})^2(1-\varepsilon^{a_j})}\deg_{C_i} x_1-\right. \\
\hfill \left. \displaystyle\frac{\varepsilon^{a_j}}{(1-\varepsilon^{-a_j})(1-\varepsilon^{a_j})^2}\deg_{C_i} x_2 \right),
\end{array}$$
where $(i,j)\in\{(1,2),\ (2,3), \ (3,1)\}$.

Since $c_1(\bn_{C_i/Z})=K_{C_i}-K_Z|_{C_i},$
we get $\deg_{C_i} x_1+\deg_{C_i} x_2=2(g_i-1)-\deg K_Z|_{C_i}.$
Then
$$\begin{array}{rcl}
\chi(X,\bl_{n}) & = &\chi(\bigo_X) + {\displaystyle\frac{1}{s}\sum_{\substack{\ep\in\mu_s \\ \ep^{a_i}\ne1\,\forall i = 1,2,3}}
\frac{\varepsilon^{-n}-1}{(1-\varepsilon^{-a_1})(1-\varepsilon^{-a_2})(1-\varepsilon^{-a_3})}}+\\
\\
&&\displaystyle \sum_{i=1}^3 \frac{\ol{nk_i}(\alpha_i-\ol{nk_i})}{4s}\deg K_Z|_{C_i} + \\
\\
&&\displaystyle\frac{1}{12s}\ol{nk_i}(\alpha_i-\ol{nk_i})(\alpha_i-2\cdot\ol{nk_i})(\deg_{C_i} x_1-\deg_{C_i} x_2),
\end{array}$$
where $\overline{\phantom{\Sigma}}$ denotes the smallest residue mod $\alpha_i$
and $k_i$ are defined by $k_ia_j=1$ mod $\alpha_i$. for $(i,j)$ as above.

We can think of $\chi(\bigo_X)$ as the RR-type expression for the pair $(X,L_{n})$, where 
$L_{n}=0\in \mbox{ Div}X\otimes\bq$ is the $\bq-$divisor corresponding to $\bl_{n}$. 
The rest of the formula is a sum of contributions coming from 
\begin{itemize}
\item dissident points of type $\frac{1}{s}(a_1,a_2,a_3),$ and
\item dissident curves of type $\frac{1}{\alpha_i}(a_j,-a_j),$
\end{itemize}
which are the singularities of~$(X,L_{n})$.

\begin{remark}\rm The attentive reader will notice that the argument in the last 
two sections is completely different: for $A_r$ curves, we use resolution of 
singularities and computations on ruled surfaces, whereas for (isolated or 
dissident) singular points, we use an auxiliary cover and equivariant 
Riemann--Roch. At least for curves without dissident points, the auxiliary cover 
construction also works, though with a twist. For an $A_r$ curve, a cyclic cover 
may not exist, but a $\mu_r\oplus\mu_r$ cover necessarily does and the argument 
goes through. The invariant $N_C$ of $X$ (the only term not explicitly defined in 
Theorem~\ref{thmdis}) depends only on the first order neighbourhood of the curve 
$C$ and can be expressed explicitly in terms of the genus C, and the singularity 
type and splitting of the normal bundle of the curve $\tilde C$ over $C$ 
in the auxiliary cover. For details, consult~\cite{anita}. 
\end{remark} 

\subsection{Conclusion of the proof} 

\begin{prfth}~\ref{thmdis}\hspace{5mm}
By Corollary~\ref{rrformm}, we only need to add the 
contributions from the dissident and isolated  
singular points to the RR formula. These were computed above. 

Let $C$ be a $_{\ol{mk}}\bigl(\frac{1}{r}(1,-1)\bigr)$ curve for $\bigo_X(mD)$ 
with dissident points
$$\{P_{\lambda}\ |\ \lambda\in\Lambda\} \mbox{ of types }\, _{\ol{mn_{\lambda}}}
\Bigl(\frac{1}{s_{\lambda} }(a_{1\lambda},a_{2\lambda},a_{3\lambda})\Bigr).$$ 
Note that for every $\lambda$ there exists $i$ such that 
$\hcf(a_{i\lambda},s_{\lambda})=r$, since $P_{\lambda}\in C$. Moreover,
$mn_{\lambda}a_{j\lambda}^{-1}=mk$ mod $r$. 

The contribution to RR from this locus is a sum of 
$$ \frac{\ol{mk}(r-\ol{mk})(r-2\cdot \ol{mk})}{6r^2}\left( (r-2)(1-g)+
\textstyle\frac{r-2}{2}K_Y\,^C\!E_1\,^C\!E_2+\,^C\!E_{r-2}^2\,^C\!E_{r-1}\right)$$
by Corollary~\ref{rrformm}, and of  
$$\begin{array}{l}
\frac{1}{12s_{\lambda}}\ol{mk}(r-\ol{mk})(r-2\cdot\ol{mk})(\deg_{C_{i\lambda}} x_1-\deg_{C_{i\lambda}} x_2)\\
\\
\mbox{}+\displaystyle\frac{1}{s_{\lambda}}\sum_{\substack{\ep\in\mu_{s\lambda} \\ \ep^{a_{i\lambda}}\ne1\,\forall i = 1,2,3}} \frac{\varepsilon^{-mn_{\lambda}}-1}{(1-\varepsilon^{-a_{1\lambda}})(1-\varepsilon^{-a_{2\lambda}})(1-\varepsilon^{-a_{3\lambda}})}.
\end{array}$$
for every point $P_{\lambda}$.
The first two rows of the sum can be written together as  
$$ \frac{\ol{mk}(r-\ol{mk})(r-2\cdot \ol{mk})}{12 r^2\tau_C}N_C,$$
where $\tau_C=\frac{1}{r}\lcm\{s_{\lambda}\}$ is the index of $C$ 
and the integer $N_C$ is an invariant of $X$ in the neighbourhood of $C$. 
\end{prfth}
\ethrm

\section{\cy threefolds}
\label{caya}

\subsection{The Hilbert series} 

When $X$ is a \cy threefold, the RR formula assumes a much more compact form.

\begin{corollary}
\label{thmcordis}
Let $X$ be a \cy threefold and let $D$ be a $\bq$-Cartier divisor.
Assume that $(X,D)$ is quasi-smooth and well formed with the 
following singularities: 
\begin{itemize}
\item points $Q$ of type $_n\bigl(\frac{1}{s}(a_1,a_2,a_3)\bigr),$
\item curves $C$ of generic type $_k\bigl(\frac{1}{r}(1,-1)\bigr)$ with index $\tau_C$. 
\end{itemize}
Then for all positive integers $m$,
$$\begin{array}{rcl}
h^0(X,\bigo_X(mD))&= & \displaystyle\frac{1}{6}m^3D^3+m\frac{D\cdot c_2(X)}{12}  \\
 & &\\
 & &\mbox{}+\sum_{Q}c_Q(mD)+\sum_{C}s_C(mD),
\end{array}$$ 
where
$$c_Q(mD)=
\frac{1}{s}\sum_{\substack{\ep\in\mu_s \\ \ep^{a_i}\ne1\,\forall i = 1,2,3}}
\frac{\varepsilon^{-nm}-1}{(1-\varepsilon^{-a_1})(1-\varepsilon^{-a_2})(1-\varepsilon^{-a_3})}$$
and
$$ s_C(mD)= -m\frac{\ol{mk}(r-\ol{mk})}{2r}\deg D|_C + \frac{\ol{mk}(r-\ol{mk})(r-2\cdot\ol{mk})}{12r^2\tau_C}N_C,$$
where $\overline{\phantom{\Sigma}}$ denotes the smallest residue mod $r$ and the 
integer $N_C$ is an invariant of $X$ in the neighbourhood of $C.$
\end{corollary}

\begin{prf}
This is a direct corollary of Theorem~\ref{thmdis}. By Kodaira vanishing, the
higher cohomologies of $\bigo_X(mD)$ vanish for $m>0$, and $\chi(\bigo_X)=0$ 
by Serre duality. 
\end{prf}
\ethrm

\begin{remark}\rm
Observe that for a $\frac{1}{2}(1,1)$ curve of singularities, 
the term involving $N_C$ vanishes. The Riemann--Roch formula for this special case 
was already proved in~\cite{balazs2}.
\end{remark}

\begin{corollary}
\label{cyrhil}
If $(X,D)$ satisfies the conditions of Corollary~\ref{thmcordis}, then the 
Hilbert series 
$$P_X(t)=1+\sum_{m=1}^{\infty}h^0(X,\bigo_X(mD))t^m.$$
can be written in the compact form
$$\begin{array}{rcl}
P_X(t)&=&
1+\displaystyle\frac{D^3}{6}\cdot \frac{t^3+4t^2+t}{(1-t)^4}+\frac{c_2(X)\cdot D}{12}\cdot \frac{t}{(1-t)^2}\\
 && \\
 & &\mbox{}+\sum_Q \tilde{P}_Q(t)+\sum_C \tilde{P}_C(t), 
\end{array}$$
where every curve $C$ contributes 
$$\begin{array}{rcl}
\tilde{P}_C(t)&=&
  -\deg D|_C\left( \displaystyle\frac{1}{1-t^r}\sum_{i=1}^{r-1}i\frac{\ol{ik}(r-\ol{ik})}{2r}t^i
+ \displaystyle\frac{rt^r}{(1-t^r)^2}\sum_{i=1}^{r-1}\frac{\ol{ik}(r-\ol{ik})}{2r}t^i \right)\\
 & &\\
 &&\mbox{}+ \displaystyle\frac{N_C}{12 r^2 \tau_c}\cdot \frac{1}{1-t^r}\sum_{i=1}^{r-1}\ol{ik}(r-\ol{ik})(r-2\cdot\ol{ik})t^i
\end{array}$$
and every singular point $Q$ contributes 
$$\tilde{P}_Q(t)=\frac{1}{1-t^s}\sum_{i=1}^{s-1}c_Q(iD)t^i.$$
Here 
$$c_Q(iD)=\frac{1}{s}\sum_{\substack{\ep\in\mu_s \\ \ep^{a_i}\ne1\,\forall i = 1,2,3}}
\frac{\varepsilon^{-ni}-1}
{(1-\varepsilon^{-a_1})(1-\varepsilon^{-a_2})(1-\varepsilon^{-a_3})}.$$ 
\end{corollary} 

\begin{prf} This follows directly from Corollary~\ref{thmcordis}, 
using elementary summation formulae of power series. 
\end{prf}
\ethrm

\subsection{Examples} 
\label{examples}

We apply Corollary~\ref{cyrhil} to the construction of new families of 
projective Calabi--Yau threefolds as follows. In the first step, we compute 
the Hilbert series $P(t)$ from Corollary~\ref{cyrhil} using the following 
input data:
\begin{itemize}
\item integers $h^0(X,D)$ and $h^0(X,2D);$
\item points $\left\{{}_n\bigl(\frac{1}{s}(a_1,a_2,a_3)\bigr)\right\};$
\item curves $\left\{{}_k\bigl(\frac{1}{r}(1,-1)\bigr)\mbox{ of degree } \deg D|_C, \mbox{ index } \tau_C \mbox{ and invariant } N_C \right\}.$
\end{itemize}
In the next step, we look for a set of weights $w_0,\ldots,w_n$ such that
$$Q(t)=P(t)\prod_{k=0}^{n}(1-t^{w_k})$$
is a polynomial. Then a plausible guess is that $(X,D)$ can be embedded in the 
ambient space $\bp(w_0,\ldots,w_n)$, and the shape of $Q(t)$ will suggest 
a set of generators and relations for the defining ideal. 
Compare~\cite{altinok, abr} for the philosophy and some explicit examples 
of this type of argument. As discussed in~\cite{abr} in great detail, 
this procedure is best done by a computer; we thank Gavin Brown for providing 
a computer program written in the programming language of the Magma computer 
algebra system~\cite{magma}. 

Lists of complete intersection \cy threefolds in weight\-ed projective spaces 
can certainly be generated in this way; such Calabi--Yau manifolds have been
listed by direct methods and extensively studied in the literature, for example 
in~\cite{cand1, cand2, cikreuzer}. Our computer search also generates more 
interesting examples in higher codimensions. 

\begin{example}
\label{cod3ex4wg} \rm Take the input data 
\begin{itemize}
\item $h^0(X,D)=3$ and $h^0(X,2D)=6$; 
\item points $\left\{{}_2\bigl(\frac{1}{3}(1,1,1)\bigr),{}_8\bigl(\frac{1}{9}(1,3,5)\bigr)\right\}$;
\item curve $\frac{1}{3}(1,2) \mbox{ with } \deg D|_C=\frac{1}{9},\ \tau_C=3,\ N_C=22$.
\end{itemize}
The Hilbert series is
\[ P(t)=\frac{Q(t)}{(1-t)^3(1-t^3)^2(1-t^5)(1-t^9)},\]
where
$$Q(t)=-t^{23}+t^{17}+t^{15}+2t^{13}+t^{11}-t^{12}-2t^{10}-t^8-t^6+1.$$
This suggests that $(X,D)$ could be realized as a codimension 3 threefold
in the weighted projective space $\bp^6(1^3,3^2,5,9)$. 
Indeed, define~$X\subset\bp^6(1^3,3^2,5,9)$ by the vanishing of 
the submaximal Pfaffians 
\begin{eqnarray*}
\label{pfeqp}
 & &\pf_1=b_7y_1+b_5 z+x_2 v,   \\
 & &\pf_2=a_7y_1+a_5 z+x_1 v,   \\
 & &\pf_3=y_2 y_1-a_5x_2+b_5x_1,  \\
 & &\pf_4=-y_2 z-a_7x_2+b_7x_1,  \\
 & &\pf_5=y_2 v-a_7b_5+a_5b_7, 
\end{eqnarray*}
of the $5\times 5$ skew-symmetric matrix 
\[ M = \left(\begin{matrix}  0 & y_2 & a_7 & a_5 & x_1 \\ 
                             & 0 & b_7 & b_5 & x_2 \\ 
                              & & 0 & v & -z \\
                               & && 0 & y_1 \\ 
                                & & & & 0
             \end{matrix}\right).\]
Here $x_1,x_2,x_3,y_1,y_2,z, v$ are weighted variables of $\bp^6$ and $a_i, b_i$ 
are general polynomials of degree~$i$. It is easy to check that $X$ is 
quasi-smooth and well formed, and that its singular locus indeed consists of an 
isolated singular point of type $\frac{1}{3}(1,1,1)$ 
and an $A_2$ curve with a dissident point of type $\frac{1}{9}(1,3,5)$.

>From this representation, one can also prove that $X$ is a general quasilinear 
section 
$$X=(5)^2 \cap (7)^2\cap \mathcal{C}\!\wgr(2,5)\subset \bp(1^3,3^2,5^3,7^2,9).$$
Here $\mathcal{C}\!\wgr(2,5)$ denotes the projective cone over $\wgr(2,5)$, 
the weighted Grassmannian as defined in \cite{wg}.
\end{example}

\begin{example}
\label{tocod4exa}\rm Take the input data 
\begin{itemize}
\item $h^0(X,D)=2$ and $h^0(X,2D)=4$;
\item point $_4\bigl(\frac{1}{5}(1,1,3)\bigr)$;
\item curve $\frac{1}{3}(1,2) \mbox{ with } \deg D|_C=1,\ \tau_C=1,\ N_C=12$.
\end{itemize}
The computer output for $P(t)$ is
\[ P(t)=\frac{t^{21}-3t^{15}-3t^{14}-3t^{13}+2t^{12}+6t^{11}+6t^{10}+2t^9-3t^8-3t^7-3t^6+1}{(1-t)^2(1-t^2)(1-t^3)^4(1-t^5)}.\]
This suggests that $(X,D)$ could be realized as a codimension 4 threefold in the weighted 
projective space $\bp^7(1^2,2,3^4,5)$, defined by nine relations of degrees $6,6,6$, $7,7,7$, 
$8,8,8$ respectively. The existence of $X$ can be proved by the Type I 
unprojection \cite{altinok, stavmile} starting from a codimension 3 variety 
$Y \subset \bp^6(1,1,2,3,3,3,3)$ defined by a set of Pfaffians as in Example~\ref{cod3ex4wg}.
We omit the details.
\end{example}

\begin{example}
\label{ex3wg} \rm Take the input data 
\begin{itemize}
\item $h^0(X,D)=2$ and $h^0(X,2D)=7$;
\item point ${}_3\bigl(\frac{1}{4}(2,3,3)\bigr)$;
\item curve $\frac{1}{2}(1,1) \mbox{ with } \deg D|_C=\frac{7}{4},\ \tau_C=2.$
\end{itemize}
The closed form of the Hilbert series is
\[ \frac{1-3t^4-4t^5-t^6+6t^7+6t^8+2t^9-2t^{11}-6t^{12}-6t^{13}+t^{14}+4t^{15}+3t^{16}-t^{20}}{(1-t)^2(1-t^2)^4(1-t^3)^2(1-t^4)}.\]
This suggests that $(X,D)$ could be realized as a codimension 5 threefold in the weighted 
projective space $\bp^8(1^2,2^4,3^2,4)$. This turns out to be true indeed, although 
the details are somewhat tedious. More conceptually, following~\cite{wg}, one can show that 
$(X,D)$ can be realized as a general quasilinear section 
$$X=\wogr(5,10) \cap(2)^2\cap (3)^4 \cap (4)\subset \bp(1^2,2^6,3^6,4^2), $$
where $\wogr(5,10)$ denotes a weighted orthogonal Grassmannian~\cite{wg}.
\end{example}

\begin{remark} \rm Two questions arise at this point: to what extent are these 
Calabi--Yau threefold families ``new'', and how many families can one construct 
this way. The second question will be discussed elsewhere; to give away the 
(negative) punchline, however hard we try, so far we have been unable 
to construct infinitely many families. In answer to the first question, 
it is certainly possible that our examples are birational to one of the 
gigantic number of toric complete intersection 
Calabi--Yau threefolds~\cite{baty, bb, torickreuzer}. However, 
our projective descriptions, in relatively simple ambient varieties 
such as weighted Grassmannians~\cite{wg} and the 
``universal Type I unprojection''~\cite{stavmile}, are new. Such descriptions 
are quite pretty in themselves, and can also be used for various purposes 
such as computing their Hodge numbers and studying mirror symmetry 
phenomena (details to follow). For these reasons alone, we believe 
they deserve some attention. 
\end{remark} 
\vspace{0.1in}

\paragraph{Acknowledgements} We thank Miles Reid for his continued attention 
and helpful advice during the course of this project, Gavin Brown for providing 
us with the Magma computer program~\cite{magma} used to find the examples 
presented in Section~\ref{examples}, and the referee for a very careful reading
of our manuscript. We also thank the Isaac Newton Institute, Cambridge 
for hospitality while part of this research was conducted, and the 
Mathematics Institute of the University of Warwick for supporting the 
research of A. B. with a Special Research Studentship.

\noindent {\small \sc Department of Mathematics, University of Ljubljana

\noindent Jadranska 19, Ljubljana 1000, Slovenia 

\noindent E-mail address: {\tt anita.buckley@fmf.uni-lj.si}

\vspace{0.1in}

\noindent {\rm and}

\vspace{0.1in}

\noindent Department of Mathematics, Utrecht University

\noindent PO. Box 80010, NL-3508 TA Utrecht, The Netherlands

\noindent E-mail address: \tt szendroi@math.uu.nl}

\end{document}

%% file: conf.pstex_t
\begin{picture}(0,0)%
\includegraphics{conf.pstex}%
\end{picture}%
\setlength{\unitlength}{4144sp}%
\begingroup\makeatletter\ifx\SetFigFont\undefined%
\gdef\SetFigFont#1#2#3#4#5{%
  \reset@font\fontsize{#1}{#2pt}%
  \fontfamily{#3}\fontseries{#4}\fontshape{#5}%
  \selectfont}%
\fi\endgroup%
\begin{picture}(4613,3036)(383,-2590)
\put(4861,-2536){\makebox(0,0)[lb]{\smash{\SetFigFont{12}{14.4}{\familydefault}{\mddefault}{\updefault}{\color[rgb]{0,0,0}$C$}%
}}}
\put(586,-2536){\makebox(0,0)[lb]{\smash{\SetFigFont{12}{14.4}{\familydefault}{\mddefault}{\updefault}{\color[rgb]{0,0,0}$E_1$}%
}}}
\put(2566,-2536){\makebox(0,0)[lb]{\smash{\SetFigFont{12}{14.4}{\familydefault}{\mddefault}{\updefault}{\color[rgb]{0,0,0}$E_{r-1}$}%
}}}
\put(1531,-2536){\makebox(0,0)[lb]{\smash{\SetFigFont{12}{14.4}{\familydefault}{\mddefault}{\updefault}{\color[rgb]{0,0,0}$\ldots$}%
}}}
\put(4996,-1366){\makebox(0,0)[lb]{\smash{\SetFigFont{12}{14.4}{\familydefault}{\mddefault}{\updefault}{\color[rgb]{0,0,0}$P$}%
}}}
\put(721,-1276){\makebox(0,0)[lb]{\smash{\SetFigFont{12}{14.4}{\familydefault}{\mddefault}{\updefault}{\color[rgb]{0,0,0}$F_1$}%
}}}
\put(2521,-1276){\makebox(0,0)[lb]{\smash{\SetFigFont{12}{14.4}{\familydefault}{\mddefault}{\updefault}{\color[rgb]{0,0,0}$F_{r-1}$}%
}}}
\end{picture}

%% file: rr.bbl
\begin{thebibliography}{99} 
\addcontentsline{toc}{section}{References}
{\small
\bibitem{altinok} S. Alt\i nok, {\it Graded rings corresponding to polarised K3 surfaces and $\bq$-Fano 3-folds}, PhD. Thesis, University of Warwick, 1998.
\bibitem{abr} S. Alt\i nok, G. Brown and M. Reid, {\it Fano 3-folds, $K3$ surfaces and graded rings}, in: Topology and geometry: commemorating SISTAG (ed. A. J. Berrick et al), 25--53, Contemp. Math., 314, Amer. Math. Soc., Providence, RI, 2002.
\bibitem{as} M. F. Atiyah and G. B. Segal, {\it The index of elliptic operators. II}, Ann. of Math. {\bf 87} (1968) 531--545.
\bibitem{atiyah} {M. F.Atiyah and I. M. Singer}, {\it The index of elliptic operators. III}, Ann. of Math. {\bf 87} (1968) 564--604.
\bibitem{baty} V. V. Batyrev, {\it Dual polyhedra and mirror symmetry for Calabi-Yau hypersurfaces in toric varieties}  J. Algebraic Geom. {\bf 3} (1994) 493--535. 
\bibitem{bb} L. A. Borisov, {\it Towards the mirror symmetry for Calabi-Yau complete intersections in Gorenstein toric Fano varieties}, alg-geom/9310001.
\bibitem{magma} W. Bosma, J. Cannon and C.~Playoust, {\it The Magma algebra system I: The user language}, J. Symb. Comp. {\bf 24} (1997) 235--265. \newline See also http://www.maths.\linebreak[2]usyd.edu.au:8000/u/magma
\bibitem{anita} A. Buckley, {\it Orbifold Riemann--Roch for threefolds and applications to Calabi--Yaus}, PhD. Thesis, University of Warwick, 2003.
\bibitem{cand1} P. Candelas, A. M. Dale, C. A. L\"utken and  R. C. Schimmrigk, {\it Complete intersection Calabi-Yau manifolds}, Frontiers of high energy physics (London, 1986, ed. I. G. Halliday), 88--134, Hilger, Bristol, 1987. 
\bibitem{cand2} P. Candelas, M. Lynker and R. Schimmrigk, {\it Calabi-Yau manifolds in weighted ${\mathbb P}^4$}, Nuclear Phys. B {\bf 341} (1990) 383--402.
\bibitem{wg}  {A. Corti} and {M. Reid}, {\it Weighted Grassmannians}, Algebraic geometry. A volume in memory of Paolo Francia (ed. M. C. Beltrametti et al), 141--163, de Gruyter, Berlin, {2002}. 
\bibitem{fletcher1} A. R. Fletcher, {\it Contributions to Riemann-Roch on projective $3$-folds with only canonical singularities and applications} Algebraic geometry, Bowdoin, 1985 (Brunswick, Maine, 1985, ed. S. J. Bloch), 221--231, Proc. Sympos. Pure Math., 46, Part 1, Amer. Math. Soc., Providence, RI, 1987.
\bibitem{fletcher} {A. R. Fletcher}, {\it Working with weighted complete intersections}, Explicit birational geometry of 3--folds (ed. A. Corti and M. Reid), 101--173, London Math. Soc. Lecture Note Ser., 281, Cambridge Univ. Press, Cambridge, {2000}. 
\bibitem{morrison} {K. Intriligator, D. R. Morrison} and {N. Seiberg}, {\it Five--dimensional supersymmetric gauge theories and degenerations of Calabi--Yau spaces}, Nuclear Phys. B {\bf 497} (1997) 56--100.
\bibitem{kawasaki} Kawasaki, T. {\it The Riemann-Roch theorem for complex $V$-manifolds}, Osaka J. Math.  {\bf 16}  (1979) 151--159.
\bibitem{torickreuzer} M. Kreuzer and H. Skarke, {\it Complete classification of reflexive polyhedra in four dimensions}, Adv. Theor. Math. Phys. {\bf 4} (2002) 1209-1230.
\bibitem{cikreuzer} {M. Kreuzer, E. Riegler} and {D. A. Sahakyan}, {\it Toric complete intersections and weighted projective space}, J. Geom. Phys. {\bf 46} (2003) 159--173.
\bibitem{stavmile} {S. Papadakis} and {M. Reid}, {\it Kustin--Miller unprojection without complexes}, J. Algebraic Geom. {\bf 13} (2004) 563--577.
\bibitem{3-folds}  {M. Reid}, {\it Canonical 3--folds}, Journ\'{e}es de g\'{e}om\'{e}trie alg\'{e}brique d'Angers (ed. A. Beauville), 1979, 273--310.
\bibitem{ypg}  {M. Reid}, {\it Young person's guide to canonical singularities},  Algebraic geometry (Bowdoin 1985, ed. S. Bloch), 345--414, Proc. Sympos. Pure Math. 46, Part 1, AMS, Providence, RI, 1987.
\bibitem{balazs2} {B. Szendr\H oi}, {\it Calabi--Yau threefolds with a curve of singularities and counter\-examples to the Torelli problem II}, Math. Proc. Cambridge Philos. Soc. {\bf 129} (2000) 193--204.
\bibitem{toen} {B. Toen}, {\it Th\'eor\'emes de Riemann-Roch pour les champs de Deligne-Mumford}, $K$-Theory {\bf 18} (1999) 33--76.
}
\end{thebibliography}
